\documentclass{article}
\newcommand{\R}{{\rm I\kern-.1567em R}}
\newtheorem{thm}{Theorem}

\begin{document}

\title{On stationary solutions to reaction-diffusion equations with nonlocal term
\footnote{\uppercase{T}his work is partially supported 
by the \uppercase{P}olish \uppercase{M}inistry of \uppercase{S}cience project 
\uppercase{N N}201 418839.}}

\author{R. STA\'NCZY \\
\small Instytut Matematyczny, Uniwersytet Wroc{\l}awski, \\
\small pl. Grunwaldzki 2/4, 50--384 Wroc{\l}aw, Poland;\\
\small{\tt stanczr@math.uni.wroc.pl}}



\maketitle

\abstract{In this paper we obtain the existence of a radial solution
for some elliptic nonlocal problem with constraints. The problem arises
from some reaction-diffusion equation modelling among others system of 
self-gravitating particles when one looks for its stationary solutions.  
}

\section{Introduction}

Consider the following parabolic--elliptic system, motivated by some astrophysical
models of self--gravitating particles (cf. \cite{CSR}) and derived by statistical mechanics approach
(cf. \cite{BS-K}),
\begin{eqnarray} \label{par}
\begin{array}{ll}
n_t = \nabla\cdot N \left(\theta P'\ \nabla n + n\ \nabla\varphi \right) \;\;\mbox{ in }\;\; \Omega\times 
(0,\infty)\,,\\
\Delta\varphi = n  \;\;\mbox{ in }\;\; \Omega\times (0,\infty)\,,\\
\left(\theta P'\ \nabla n + n\ \nabla\varphi \right)\cdot\bar{\nu} = \varphi = 0  \;\;\mbox{ on }\;\;  
\partial\Omega\times (0,\infty)\,,\\
n(0) = n_0 \ge 0 \;\;\mbox{ in }\;\; \Omega\subset {\R}^d\,, 
\end{array}
\end{eqnarray}
where $n=n(x,t)$ is a nonnegative density of the particles inducing
the gravitational potential $\varphi = \varphi(x,t)$, the pressure term $P$ 
depending on $n\theta^{-d/2}$ and any positive (otherwise irrelevant for stationary problem) 
coefficient $N$. Note that total mass of the system 
is conserved, i.e., $\int_{\Omega} n(x,t) \, dx =m$.
The system can be reduced to one nonlocal equation of reaction-diffusion type
if we plug $\varphi = -\Delta^{-1}n$ into the first equation.

Stationary solutions can be derived by multiplication of the first eqution in 
(\ref{par}) by $\theta H+\varphi$ (where $H'(z)=P'(z)/z$) and integration to yield
\begin{eqnarray}\label{ell}
\begin{array}{ll}
\Delta \varphi = n = H^{-1}(c-\varphi),\\
\int_{\Omega} n = m,
\end{array}
\end{eqnarray}

For radially symmetric $\Omega,$ e.g. a unit ball of radius one $B(0,1)$ 
looking for radial solutions in the form of integrated densities $Q(r)=\int_{B(r)} n(x) dx$
we are reduced to the following singular boundary value problem (cf. for $R=Id$ \cite{BDEMN}, \cite{BN})

\begin{eqnarray}\label{rad}
\begin{array}{ll}
-Q''(r)+(d-1)r^{-1}Q'(r)=R(Q'(r)r^{1-d}\sigma_d^{-1})Q (r){\rm \;\; for \;\;} r\in (0,1),\\
Q(0)=0,  {\rm \;\; and \;\;} Q(1)=m,
\end{array}
\end{eqnarray}
for the given mass parameter $m>0$ and $\sigma_d$ being the measure of the unit
sphere in ${\R}^d,$ $\theta=1$ and $R(z)=1/H'(z)$. Indeed, differentiating
$H(n)=c-\varphi$ and using radial symmetry yields the claim.

Next integrating the above equation twice and using the boundary conditions 
we can reduce the above problem to looking for fixed points of the following
operator
\begin{equation}
{\mathcal T}Q(r)=mr^d+\frac{1}{d}\int_0^1 R(Q's^{1-d}\sigma_d^{-1})s^{1-d}Q(s)G(r,s)\,ds,
\end{equation}
where a symmetric function $G$ is given by
\begin{equation}
G(r,s)=r^d(1-s^d)  {\rm \;\; for \;\;} s>r.
\end{equation}

\section{Main results}

Now, we are ready to formulate the main result of this paper.

\begin{thm}
Assume that the function $R$ is globally Lipschitz continuous with constant $L$ and $R(0)=0$.
Then, for sufficiently small positive mass parameter $m$ there exists at least one solution
to $(\ref{rad}).$ Moreover, if $R=Id$ then the function $Q(r)r^{2-d}$ is non-decreasing.
\end{thm}

\noindent
{\bf Proof.}
Next note that if $R$ is locally Lipschitz continuous then in the weighted sup norm 
$|Q|_\alpha=\sup_{r\in (0,1)}|Q(r)r^{\alpha}|$ for $\alpha\le 0$ the estimates hold
\begin{eqnarray*}
\begin{array}{ll}
|{\mathcal T}Q|_{2-d}\le A_1 |Q|_{2-d}|Q'|_{3-d}+m,\\
|\left({\mathcal T}Q\right)'|_{3-d}\le A_2 |Q|_{2-d}|Q'|_{3-d}+md,
\end{array}
\end{eqnarray*}
for some constants $A_1,A_2$ depending on the Lipschitz constant of $R$ (call it $L$) and the dimension
of the space ${\R}^d.$ To be more specific $2A_1\sigma_d(d-2)=L$ and $2A_2\sigma_d(d-2)=L(d+4)$.
Moreover, similarly one can get the estimates 
\begin{eqnarray*}
\begin{array}{ll}
|{\mathcal T}Q-{\mathcal T}S|_{2-d} \le  A_3 \max\{ |Q-S|_{2-d}|Q'|_{3-d}, |Q'-S'|_{3-d}||S|_{2-d}\},\\
|\left({\mathcal T}Q\right)'-\left({\mathcal T}S\right)'|_{3-d} \le A_4  \max\{ |Q-S|_{2-d}|Q'|_{3-d}, |Q'-S'|_{3-d}||S|_{2-d}\}
\end{array}
\end{eqnarray*}
where the constants are defined as: $2A_3\sigma_d(d-2)=L+1$ and $2A_4\sigma_d(d-2)=L(d+4)+1$.
Thus ${\mathcal T}$ is a contraction on some ball $B(0,\rho)$ in the following subspace of $C^1$ -
$C^1_d=\{Q\in C^1:|Q|_{2-d}<\infty, |Q'|_{3-d}<\infty \}$ with the weighted sup norm 
$$
\|Q\|=\max\{|Q|_{2-d},|Q'|_{3-d}\}
$$
provided that 
\begin{equation}\label{sma}
\rho(L(d+4)+1)<2\sigma_d(d-2), {\rm \;\;} \rho^2\frac{L(d+4)}{2\sigma_d(d-2)}+md\le \rho,
\end{equation}
yielding the existence of a fixed point for ${\mathcal T}$ and thus a solution to $(\ref{rad})$.
Note that if $m$ is sufficiently small then one can find $\rho \in (\rho^1_m,\rho^2_m)$ satisfying (\ref{sma})
for some suitably chosen $\rho^1_m, \rho^2_m$.

To prove another part of the claim for $R=Id$ we use the property of invariance of
the cone $\{Q\in C^1_d:\left(Q(r)r^{2-d}\right)'\ge 0 \}$ under
the action of the operator ${\mathcal T}$ provided mass is sufficiently small
($m<2\sigma_d$).


\end{document}